\documentclass{amsart}

\usepackage{amssymb, amsthm}

\newtheorem{thm}{Theorem}[section]
\newtheorem{lem}[thm]{Lemma}

\theoremstyle{definition}
\newtheorem{defn}[thm]{Definition}

\theoremstyle{remark}
\newtheorem{rem}[thm]{Remark}

\def\zero{^{\circ}}

\begin{document}

\title{A Computation of tight closure in diagonal hypersurfaces}

\author{Anurag K. Singh}

\address{Department Of Mathematics, The University Of Michigan, East Hall,
525 East \linebreak University Avenue, Ann Arbor, MI 48109-1109}

\email{binku@math.lsa.umich.edu}

\date{October 16, 1997}
\maketitle

\section{Introduction}

The aim of this paper is to settle a question about the tight closure of the
ideal $(x^2, y^2, z^2)$ in the ring $R=K[X,Y,Z]/(X^3+Y^3+Z^3)$ where $K$ is a 
field of prime characteristic $p \ne 3$. (Lower case letters denote the images 
of the corresponding variables.) M.~McDermott has studied the tight closure of
various irreducible ideals in $R$, and has established that 
$xyz \in (x^2, y^2, z^2)^*$ when $p < 200$, see \cite{moira}. The general case 
however existed as a classic example of the difficulty involved in tight closure 
computations, see also \cite[Example 1.2]{huneke}. We show that 
$xyz \in (x^2, y^2, z^2)^*$ in arbitrary prime characteristic $p$, and 
furthermore establish that $xyz \in (x^2, y^2, z^2)^F$ whenever $R$ is not 
F--pure, i.e., when $p \equiv 2 \mod 3$.  We move on to generalize these results 
to the diagonal hypersurfaces $R=K[X_1, \dots, X_n]/(X_1^n + \dots + X_n^n)$. 

\medskip

These issues relate to the question whether the tight closure $I^*$ of an ideal 
$I$ agrees with its plus closure, $I^+ = IR^+ \cap R$, where $R$ is a
domain over a field of characteristic $p$ and $R^+$ is the integral closure of
$R$ in an algebraic closure of its fraction field. In this setting, we may
think of the Frobenius closure of $I$ as $I^F = IR^{\infty} \cap R$ where 
$R^{\infty}$ is the extension of $R$ obtained by adjoining $p^e\,$th roots of
all nonzero elements of $R$ for $e \in \mathbb N$. It is not difficult to see
that $I^+ \subseteq I^*$, and equality in general is a formidable open question. 
It should be mentioned that in the case when $I$ is an ideal generated by part
of a system of parameters, the equality is a result of K.~Smith, see 
\cite{Sminv}. In the above ring $R=K[X,Y,Z]/(X^3+Y^3+Z^3)$ where $K$ is a field 
of characteristic $p \equiv 2 \mod 3$, if one could show that $I^* = I^F$ for 
an ideal $I$, a consequence of this would be  
$I^F \subseteq I^+ \subseteq I^* = I^F$, by which $I^+ = I^*$. McDermott does 
show that $I^* = I^F$ for large families of irreducible ideals and our result 
$xyz \in (x^2, y^2, z^2)^F$, we believe, fills in an interesting remaining case. 

\section{Definitions}
Our main reference for the theory of tight closure is \cite{HHjams}. We next
recall some basic definitions.

\medskip

Let $R$ be a Noetherian ring of characteristic $p > 0$. We shall always use the 
letter $e$ to denote a variable nonnegative integer, and $q$ to denote the
$e\,$th power of $p$, i.e., $q=p^e$. We shall denote by $F$, the Frobenius 
endomorphism of $R$, and by $F^e$, its $e\,$th iteration, i.e., $F^e(r)=r^q$. 
For an ideal $I=(x_1, \dots, x_n) \subseteq R$, we let 
$I^{[q]}=(x_1^q, \dots, x_n^q)$. Note that $F^e(I)R=  I^{[q]}$, where $q=p^e$, 
as always. 

\medskip

We shall denote by $R\zero$ the complement of the union of the minimal 
primes of $R$.

\begin{defn} 
A ring $R$ is said to be {\it F--pure }\/ if the Frobenius homomorphism
$F: M \to M\otimes_R F(R)$ is injective for all $R$--modules $M$.

\medskip

For an element $x$ of $R$ and an ideal $I$, we say that $x \in I^F$, the 
{\it Frobenius closure }\/ of $I$, if there exists $q=p^e$ such that
$x^q \in I^{[q]}$. A normal domain $R$ is F--pure if and only if
for all ideals $I$ of $R$, we have $I^F = I$.

\medskip

We say that $x \in I^*$, the {\it tight closure }\/ of $I$, if there exists 
$c \in R\zero$ such that $cx^q \in I^{[q]}$ for all $q=p^e \gg 0$. 
\end{defn}

It is easily verified that $I \subseteq I^F \subseteq I^*$.  Furthermore, $I^*$ 
is always contained in the integral closure of $I$ and is frequently much 
smaller.

\section{Preliminary computations}

We record some determinant computations we shall find useful. Note that for
integers $n$ and $m$ where $m \ge 1$, we shall use the notation:
$$
\binom{n}{m} = \frac{(n)(n-1) \dotsm (n-m+1)}{(m)(m-1) \dotsm (1)}.
$$

\begin{lem}
$$
\det\begin{vmatrix}
\binom{n}{a+k} & \binom{n}{a+k+1} & \hdots &  \binom{n}{a+2k} \\
\binom{n}{a+k-1} & \binom{n}{a+k} & \hdots &  \binom{n}{a+2k-1} \\
\hdotsfor4 \\
\hdotsfor4 \\
\binom{n}{a}   & \binom{n}{a+1} & \hdots &    \binom{n}{a+k} \\
\end{vmatrix} 
=\frac{\binom{n}{a+k}\binom{n+1}{a+k}\dotsm\binom{n+k}{a+k}}
{\binom{a+k}{a+k} \binom{a+k+1}{a+k}\cdots\binom{a+2k}{a+k}}.
$$
\label{det1}
\end{lem}

\begin{proof}
This is evaluated in \cite[page 682]{muir} as well as \cite{roberts}.
\end{proof}

\begin{lem}
Let $F(n,a,k)$ denote the determinant of the matrix 
$$
M(n,a,k) = 
\begin{pmatrix}
\binom{n}{a} & \binom{n}{a+1} & \hdots &  \binom{n}{a+k} \\
\binom{n+2}{a+1} & \binom{n+2}{a+2} & \hdots &  \binom{n+2}{a+k+1} \\
\hdotsfor4 \\
\hdotsfor4 \\
\binom{n+2k}{a+k} & \binom{n+2k}{a+k+1} & \hdots &  \binom{n+2k}{a+2k} 
\end{pmatrix}. 
$$
Then for $k \ge 1$ we have
$$
\frac{F(n,a,k)}{F(n+2,a+2,k-1)} = \binom{n}{a} \prod_{s=1}^{k}\prod_{r=1}^{k} 
\frac{s(s+2a-n)}{(a+r)(n-a+r)}.  
$$
Hence 
$$
F(n,a,k)=\frac{\binom{n}{a}  \binom{n+2}{a+2}   \dotsm \binom{n+2k}{a+2k}}
              {\binom{a+k}{k}\binom{a+k+1}{k-1} \dotsm \binom{a+2k-1}{1}}
 .\frac{\binom{2a-n+k}{k} \binom{2a-n+k+1}{k-1} \dotsm \binom{2a-n+2k-1}{1}}
       {\binom{n-a+k}{k}  \binom{n-a+k-1}{k-1}  \dotsm \binom{n-a+1}{1}}.  
$$
\label{det2}
\end{lem}

\begin{proof}
We shall perform row operations on $M(n,a,k)$ in order to get zero entries in 
the first column from the second row onwards, starting with the last row and 
moving up. More precisely, from the $(r+1)\,$th row, subtract the $r\,$th row 
multiplied by ${\binom{n+2r}{a+r}}/{\binom{n+2r-2}{a+r-1}}$
starting with $r=k$, and continuing until $r=2$. 
The $(r+1, s+1)\,$th entry of the new matrix, for $r \ge 1$,  is
$$
\binom{n+2r}{a+r+s} - \frac{\binom{n+2r}{a+r}}{\binom{n+2r-2}{a+r-1}}
\binom{n+2r-2}{a+r-1+s}
= \frac{s(s+2a-n)}{(a+r)(n-a+r)}\binom{n+2r}{a+r+s}. 
$$
We have only one nonzero entry in the first column, namely $\binom{n}{a}$
and so we examine the matrix obtained by deleting the first row and column.
Factoring out $s(s+2a-n)$ from each column for $s=1, \dots, k$ and 
$1/(a+r)(n-a+r)$ from each row for $r=1, \dots, k$, we see that 
$$
\text{ det }M(n,a,k) =  \binom{n}{a} \prod_{s=1}^{k}\prod_{r=1}^{k} 
\frac{s(s+2a-n)}{(a+r)(n-a+r)} \text{ det }M(n+2,a+2,k-1).
$$
The required result immediately follows. 
\end{proof}

\begin{lem}
Consider the polynomial ring $T =K[A_1, \dots, A_m]$ where $I_{r,i}$ denotes the 
ideal $I_{r,i} = (A_1^i, \dots, A_r^i)T$ for $r \le m$. Then 
$$
(A_1 \dotsm A_{r-1})^{\alpha}(A_1 +\dots +A_{r-1})^{\beta} \in 
I_{r-1, \alpha+\gamma}+(A_1+ \dots +A_{r-1})^{\alpha + \gamma}T 
$$
for positive integers $\alpha, \beta$, and $\gamma$ implies
$$
(A_1 \dotsm A_r)^{\alpha}(A_1 +\dots +A_r)^{\beta+\gamma-1} \in 
I_{r, \alpha+\gamma}+(A_1+ \dots +A_r)^{\alpha + \gamma}T. 
$$
\label{induce}
\end{lem}

\begin{proof}
Consider the binomial expansion of 
$(A_1 \dotsm A_r)^{\alpha}(A_1 +\dots +A_r)^{\beta+\gamma-1}$ into terms of the
form 
$(A_1 \dotsm A_{r-1})^{\alpha}(A_1 +\dots +A_{r-1})^{\beta+\gamma-1-j}
A_r^{\alpha+j}$. Such an element is clearly in 
$I_{r, \alpha+\gamma}$ whenever $j \ge \gamma$, and so assume $\gamma > j$. 
Now
\begin{align*}
(A_1 \dotsm A_{r-1})^{\alpha}A_r^{\alpha+j}
&(A_1 +\dots +A_{r-1})^{\beta+\gamma-1-j} \\
& \in I_{r,\alpha+\gamma}+A_r^{\alpha+j}(A_1+\dots +A_{r-1})^{\alpha+2\gamma
-1-j}T \\
& \subseteq I_{r,\alpha+\gamma}+ (A_1 +\dots +A_{r-1},A_r)^{2\alpha+2\gamma-1}T
\\
& \subseteq I_{r, \alpha+\gamma}+(A_1 +\dots +A_r)^{\alpha+\gamma}T.
\end{align*}
\end{proof}

\section{Tight closure}
We now prove the main theorem. 

\begin{thm}
Let $R=K[X_1, \dots, X_n]/(X_1^n + \dots + X_n^n)$ where $n \ge 3$ and $K$ is a 
field of prime characteristic $p$ where $p \nmid n$. Then
$$
(x_1 \dotsm x_n)^{n-2} \in (x_1^{n-1}, \dots, x_n^{n-1})^*. 
$$
\end{thm}

Note that there are infinitely many $e \in \mathbb N$ such that 
$p^e=q \equiv 1 \mod n$. 
By \cite[Lemma 8.16]{HHjams}, it suffices to work with powers of $p$ 
of this form, and show that for all such $q$ we have 
$$
(x_1 \dotsm x_n)^{(n-2)q+1} \in (x_1^{(n-1)q}, \dots, x_n^{(n-1)q}). 
$$
Letting $q=nk+1$, it suffices to show
$$
(x_1 \dotsm x_n)^{(n-2)nk} \in (x_1^{(n-1)nk}, \dots, x_n^{(n-1)nk}). 
$$
Let $A_1 = x_1^n, \dots, \ A_n = x_n^n$ and note that $A_1 + \dots + A_n=0$. 
In this notation, we aim to show 
$$
(A_1 \dotsm A_n)^{(n-2)k} \in (A_1^{(n-1)k}, \dots, A_n^{(n-1)k}). 
$$
Our task is then effectively reduced to working in the polynomial ring
$$
K[A_1, \dots, A_{n-1}] \cong K[A_1, \dots, A_n]/(A_1 + \dots + A_n)
$$
where we need to show 
$$
(A_1 \dotsm A_{n-1}(A_1+\dots +A_{n-1}))^{(n-2)k} 
\in I_{n-1,(n-1)k} + (A_1+\dots +A_{n-1})^{(n-1)k}.
$$
By repeated use of Lemma \ref{induce}, it suffices to show
$$
(A_1A_2)^{(n-2)k}(A_1+A_2)^{k} 
\in (A_1^{(n-1)k}, A_2^{(n-1)k}, (A_1+A_2)^{(n-1)k}).
$$
We have now reduced our problem to a statement about a polynomial ring in two
variables. The required result follows from the next lemma. 

\begin{lem}
Let $K[A,B]$ be a polynomial ring over a field $K$ of characteristic $p >0$
and $e$ be a positive integer such that $q=p^e \equiv 1 \mod n$. If
$q=nk+1$, we have
$$
(A,B)^{(2n-3)k} \subseteq I = (A^{(n-1)k}, \  B^{(n-1)k}, \  (A+B)^{(n-1)k}).
$$
In particular, $(AB)^{(n-2)k}(A+B)^k \in I$.
\label{general}
\end{lem}

\begin{proof}
Note that $I$ contains the following elements: $(A+B)^{(n-1)k}A^kB^{(n-3)k}$, 
\newline
$(A+B)^{(n-1)k}A^{k-1}B^{(n-3)k+1}, \dots, \ (A+B)^{(n-1)k}B^{(n-2)k}$.
We take the binomial expansions of these elements and consider them modulo
the ideal $(A^{(n-1)k}, B^{(n-1)k})$. This shows that the following elements 
are in $I:$
$$
\begin{array}{rrrrr}
\binom{(n-1)k}{k}A^{(n-1)k}B^{(n-2)k} 
&+ &\cdots
&+ &\binom{(n-1)k}{2k}A^{(n-2)k}B^{(n-1)k}, \\
\binom{(n-1)k}{k-1}A^{(n-1)k}B^{(n-2)k} 
&+ &\cdots
&+ &\binom{(n-1)k}{2k-1}A^{(n-2)k}B^{(n-1)k}, \\
\hdotsfor5 \\
\hdotsfor5 \\
\binom{(n-1)k}{0}A^{(n-1)k}B^{(n-2)k} 
&+ &\cdots 
&+ &\binom{(n-1)k}{k}A^{(n-2)k}B^{(n-1)k}.
\end{array}
$$
The coefficients of 
$A^{(n-1)k}B^{(n-2)k}, A^{(n-1)k-1}B^{(n-2)k+1}, \dots,  A^{(n-2)k}B^{(n-1)k}$  
form the matrix: 
$$
\begin{pmatrix}
\binom{(n-1)k}{k}   & \binom{(n-1)k}{k+1} &  \hdots &  \binom{(n-1)k}{2k} \\
\binom{(n-1)k}{k-1}   & \binom{(n-1)k}{k} &  \hdots &  \binom{(n-1)k}{2k-1} \\
\hdotsfor4 \\
\hdotsfor4 \\
\binom{(n-1)k}{0}   & \binom{(n-1)k}{1}   &  \hdots &  \binom{(n-1)k}{k}.
\end{pmatrix}.
$$
To show that all monomials of degree $(2n-3)k$ in $A$ and $B$ are in $I$, it 
suffices to show that this matrix is invertible. Since  $q=nk+1$ we have 
$\binom{(n-1)k+r}{k}=(-1)^k \binom{2k-r}{k}$ for $0 \le r \le k$, 
and so by Lemma \ref{det1}, the determinant of this matrix is  
$$ 
\frac{\binom{(n-1)k}{k}\binom{(n-1)k+1}{k}\cdots\binom{nk}{k}}
{\binom{k}{k} \binom{k+1}{k}\cdots\binom{2k}{k}}
= (-1)^{k(k+1)} \frac{\binom{2k}{k}\binom{2k-1}{k}\cdots\binom{k}{k}}
{\binom{k}{k} \binom{k+1}{k}\cdots\binom{2k}{k}} = 1.
$$
\end{proof}

With this we complete the proof that 
$(x_1 \dotsm x_n)^{n-2} \in (x_1^{n-1}, \dots, x_n^{n-1})^*$.

\section{Frobenius closure}

Let $R=K[X_1, \dots, X_n]/(X_1^n + \dots + X_n^n)$ as before, where the
characteristic of $K$ is $p \nmid n$. 

\begin{lem}
Let $R=K[X_1, \dots, X_n]/(X_1^n + \dots + X_n^n)$ where $K$ is a field of
characteristic $p$. Then $R$ is F--pure if and only if $p \equiv 1 \mod n$. 
\label{L:fpure}
\end{lem}

\begin{proof}
This is Proposition 5.21 (c) of \cite{HRcohom}.
\end{proof}

The main result of this section is the following theorem.

\begin{thm}
Let $R=K[X_1, \dots, X_n]/(X_1^n + \dots + X_n^n)$ where $K$ is a field of
characteristic $p$.  Then 
$$
(x_1 \dotsm x_n)^{n-2} \in (x_1^{n-1}, \dots, x_n^{n-1})^F. 
$$
if and only if $p \not\equiv 1 \mod n$. 
\end{thm}

One implication follows from Lemma \ref{L:fpure}, and so we need to consider
the case $p \not\equiv 1 \mod n$.

The case $n=3$ seems to be the most difficult, and we handle that first. 
Let $R=K[X,Y,Z]/(X^3+Y^3+Z^3)$ where $p \equiv 2 \mod 3$. We need to show that 
$xyz \in (x^2, y^2, z^2)^F$. 

\medskip

Let $A=y^3, B=z^3$ and so $A+B = -x^3$. We first show that when $p =2$, we have 
$xyz \in (x^2, y^2, z^2)^F$ by establishing that 
$(xyz)^8 \in (x^2, y^2, z^2)^{[8]}$. Note that it suffices to show that 
$(xyz)^6 \in (x^{15}, y^{15}, z^{15})$, or in other words that 
$(AB(A+B))^2 \in (A^5, B^5, (A+B)^5)$, but this is easily seen to be true.

\medskip

We may now assume $p = 6m+5$ where $m \ge 0$. We shall show that in this case 
$(xyz)^p \in (x^2, y^2, z^2)^{[p]}$, i.e., that
$$
(xyz)^{6m+5} \in (x^{12m+10}, \ y^{12m+10}, \ z^{12m+10}).
$$
Note that to establish this, it suffices to show 
$$
(xyz)^{6m+3} \in (x^{12m+9}, \ y^{12m+9}, \ z^{12m+9}), 
$$
i.e., that 
$\ (AB(A+B))^{2m+1} \in (A^{4m+3}, B^{4m+3}, (A+B)^{4m+3})$.

\medskip
\begin{lem}
Let $K[A,B]$ be a polynomial ring over a field $K$ of characteristic $p =6m+5$ 
where $m \ge 0$. Then we have
$$
(AB(A+B))^{2m+1} \in I = (A^{4m+3}, \ B^{4m+3}, \ (A+B)^{4m+3}).
$$
\label{nfp}
\end{lem}

\begin{proof}
To show that $(AB(A+B))^{2m+1} \in I$, we shall show that the following terms
grouped together symmetrically from its binomial expansion,
$$
\begin{array}{l}
f_1 = (AB)^{3m+1}(A+B), \ f_3 = (AB)^{3m}(A^3+B^3), \ \dots, \\
f_{2m+1} = (AB)^{2m+1}(A^{2m+1}+B^{2m+1}),
\end{array}
$$
are all in the ideal $I$. 
Note that $I$ contains the elements $(AB)^m(A+B)^{4m+3}$, 
$(AB)^{m-1}(A+B)^{4m+5},  \dots, (AB)(A+B)^{6m+1}, (A+B)^{6m+3}$.
We consider the binomial expansions of these elements modulo 
$(A^{4m+3}, B^{4m+3})$, and get the following elements in $I:$ 
$$
\begin{array}{rrrrrrr}
\binom{4m+3}{2m+2}f_1 &+ & \binom{4m+3}{2m+3}f_3 &+ & \dots &+ & \binom{4m+3}{3m+2}f_{2m+1}, \\
\binom{4m+5}{2m+3}f_1 &+ & \binom{4m+5}{2m+4}f_3 &+ & \dots &+ & \binom{4m+5}{3m+3}f_{2m+1}, \\
\hdotsfor7 \\
\hdotsfor7 \\
\binom{6m+3}{3m+2}f_1 &+ & \binom{6m+3}{3m+3}f_3 &+ & \dots &+ & \binom{6m+3}{4m+2}f_{2m+1}. 
\end{array}
$$
The coefficients of $f_1, f_3, \dots, f_{2m+1}$ arising from these 
terms form the matrix: 
$$
\begin{pmatrix}
\binom{4m+3}{2m+2} & \binom{4m+3}{2m+3} & \hdots &  \binom{4m+3}{3m+2} \\
\binom{4m+5}{2m+3} & \binom{4m+5}{2m+4} & \hdots &  \binom{4m+5}{3m+3} \\
\hdotsfor4 \\
\hdotsfor4 \\
\binom{6m+3}{3m+2} & \binom{6m+3}{3m+3} & \hdots &  \binom{6m+3}{4m+2} 
\end{pmatrix}.
$$
We need to show that this matrix is invertible, but in the notation of 
Lemma \ref{det2}, its determinant is $F(4m+3, 2m+2, m)$ and is easily seen 
to be nonzero. 
\end{proof}

\medskip

The above lemma completes the case $n=3$. We may now assume $n \ge 4$ and 
$p = nk+\delta$ for $2 \le \delta \le n-1$. If $k=0$ i.e., $2 \le p \le n-1$, 
we have
\begin{align*}
(x_1 \dotsm x_n)^{(n-2)p} 
& = -(x_1 \dotsm x_{n-1})^{(n-2)p}x_n^{(n-2)p-n}(x_1^n + \dots + x_{n-1}^n) \\
& \in (x_1^{(n-1)p}, \dots, x_{n-1}^{(n-1)p}).
\end{align*}

\medskip

In the remaining case, we have  $n \ge 4$ and $k \ge 1$. To prove that 
$(x_1 \dotsm x_n)^{n-2} \in (x_1^{n-1}, \dots, x_n^{n-1})^F$, we shall show
$$
(x_1 \dotsm x_n)^{(n-2)p} \in (x_1^{(n-1)p}, \dots, x_n^{(n-1)p}). 
$$
This would follow if we could show
$$
(x_1 \dotsm x_n)^{(n-2)nk} \in (x_1^{(n-1)nk+n}, \dots, x_n^{(n-1)nk+n}). 
$$
As before, let $A_1 = x_1^n, \dots, A_n = x_n^n$. It suffices to show that 
$$
(A_1 \dotsm A_n)^{(n-2)k} \in (A_1^{(n-1)k+1}, \dots, A_n^{(n-1)k+1}). 
$$
By Lemma \ref{induce}, this reduces to showing
$$
(A_1A_2)^{(n-2)k}(A_1+A_2)^k \in 
I = (A_1^{(n-1)k+1}, \ A_2^{(n-1)k+1}, \  (A_1+A_2)^{(n-1)k+1}).
$$
The only remaining ingredient is the following lemma. 

\begin{lem}
Let $K[A,B]$ be a polynomial ring over a field $K$ of characteristic $p >0$
where $p=nk+\delta$ where $n \ge 4$, $k \ge 1$ and $2 \le \delta \le n-1$. 
Then 
$$
(A,B)^{(2n-3)k} \subseteq I =(A^{(n-1)k+1}, \ B^{(n-1)k+1}, \ (A+B)^{(n-1)k+1}).
$$
In particular, $(AB)^{(n-2)k}(A+B)^k \in I$.
\end{lem}

\begin{proof} Note that $I$ contains the elements: 
$\ (A+B)^{(n-1)k+1}A^kB^{(n-3)k-1}$,
$$
(A+B)^{(n-1)k+1}A^{k-1}B^{(n-3)k}, \dots, \ (A+B)^{(n-1)k+1}B^{(n-2)k-1}. 
$$
We take the binomial expansions of these elements and consider them modulo the 
ideal $(A^{(n-1)k+1}, B^{(n-1)k+1})$. This shows that the following elements 
are in $I:$
$$
\begin{array}{rrrrr}
\binom{(n-1)k+1}{k+1}A^{(n-1)k}B^{(n-2)k} 
&+ &\cdots
&+ &\binom{(n-1)k+1}{2k+1}A^{(n-2)k}B^{(n-1)k}, \\
\binom{(n-1)k+1}{k}A^{(n-1)k}B^{(n-2)k} 
&+ &\cdots
&+ &\binom{(n-1)k+1}{2k}A^{(n-2)k}B^{(n-1)k}, \\
\hdotsfor5\\
\hdotsfor5\\
\binom{(n-1)k+1}{1}A^{(n-1)k}B^{(n-2)k} 
&+ &\cdots 
&+ &\binom{(n-1)k+1}{k+1}A^{(n-2)k}B^{(n-1)k}.
\end{array}
$$
The coefficients of 
$A^{(n-1)k}B^{(n-2)k}, A^{(n-1)k-1}B^{(n-2)k+1}, \dots,  A^{(n-2)k}B^{(n-1)k}$  
form the matrix: 
$$
\begin{pmatrix}
\binom{(n-1)k+1}{k+1} & \binom{(n-1)k+1}{k+2} &  \hdots &  \binom{(n-1)k+1}{2k+1} \\
\binom{(n-1)k+1}{k} & \binom{(n-1)k+1}{k+1} &  \hdots &  \binom{(n-1)k+1}{2k} \\
\hdotsfor4 \\
\hdotsfor4 \\
\binom{(n-1)k+1}{1}   & \binom{(n-1)k+1}{2}   &  \hdots &  \binom{(n-1)k+1}{k+1}  
\end{pmatrix}.
$$
To show that all monomials of degree $(2n-3)k$ in $A$ and $B$ are in $I$, it 
suffices to show that this matrix is invertible. The determinant of this matrix 
is  
$$ 
\frac{\binom{(n-1)k+1}{k+1}\binom{(n-1)k+2}{k+1}\cdots\binom{nk+1}{k+1}}
{\binom{k+1}{k+1} \binom{k+2}{k+1}\cdots\binom{2k+1}{k+1}}
$$
which is easily seen to be nonzero since the characteristic of the field is
$p=nk+\delta$ where $2 \le \delta \le n-1$. 
\end{proof}

\begin{rem}
It is worth noting that $xyz \in (x^2, y^2, z^2)^*$ in the ring 
$$
R=K[X,Y,Z]/(X^3+Y^3+Z^3)
$$
is, in a certain sense, unexplained. Under mild  hypotheses on a ring, tight
closure has a \lq\lq colon--capturing\rq\rq \  property: for $x_1, \dots x_n$
part of a system of parameters for an excellent local (or graded)
equidimensional ring $A$, we have  $(x_1, \dots x_{n-1}) :_A x_n \subseteq
(x_1, \dots x_{n-1})^*$ and various instances of elements being in the tight
closure of ideals are easily seen to  arise from  this colon--capturing
property.

\medskip

To illustrate our point, we recall from \cite[Example 5.7]{Ho} how 
$z^2 \in (x,y)^*$ in the ring $R$ above is seen to arise from colon--capturing. 
Consider the Segre product $T = R \# S$ where $S=K[U,V]$. Then the elements
$xv-yu$, $xu$ and $yv$ form a system of parameters for the ring $T$. This ring
is not Cohen--Macaulay as seen from the relation on the parameters:
$$
(zu)(zv) (xv-yu) = (zv)^2 (xu) - (zu)^2(yv).
$$
The colon--capturing property of tight closure shows
$$
(zu)(zv) \in \quad (xu, yv):_T(xv-yu) \quad \subseteq (xu, yv)^*.
$$
There is a retraction $R \otimes_K S \to R$ under which $U \mapsto 1$ and 
$V \mapsto 1$. This gives us a retraction from $T \to R$ which, when applied to
$(zu)(zv) \in (xu, yv)^*$, shows $z^2 \in (x,y)^*$ in $R$. 

\end{rem}

\end{document}